\documentclass[a4paper,12pt]{article}
\usepackage{amsmath,amssymb,amsfonts}
\usepackage{colordvi}
\usepackage[all]{xy}

\def\zero{\{0\}}
\def\un{{\bf 1}}
\def\dom{\backslash}
\def\mpoint{\;.}
\def\mvirg{\;,}
\def\sur{\overline}

\def\mpn{\medskip\par\noindent}

\def\mmpn{\vskip 1em minus 1em\par\noindent}

\def\smp{\smallskip\par}

\def\Br{{\rm Br}}

\def\Res{\hbox{\rm Res}}
\def\Ind{\hbox{\rm Ind}}
\def\Hom{\hbox{\rm Hom}}

\def\Inf{\hbox{\rm Inf}}

\def\Iso{\hbox{\rm Iso}}

\def\tr{\hbox{\rm tr}}
\newcommand{\gen}[1]{{<}#1{>}}

\def\normal{\mathop{\underline\triangleleft}}

\def\Z{\mathbb{Z}}

\def\Q{\mathbb{Q}}

\def\C{\mathbb{C}}

\newcommand{\sumb}[2]{\sum_{{\scriptstyle #1}\atop {\scriptstyle #2}}}
\newcommand{\sumc}[3]{\sum_{{\scriptstyle #1}\atop {{\scriptstyle #2}\atop {\scriptstyle #3}}}}

\def\pf{\par\bigskip\noindent{\bf Proof~: }}
\def\endpf{\nolinebreak~\leaders\hbox to 1em{\hss\
\hss}\hfill~\raisebox{.5ex}
{\framebox[1ex]{}}\par\bigskip}

\makeatletter
\renewenvironment{enumerate}{\ifnum \@enumdepth >3 \@toodeep\else
       \advance\@enumdepth \@ne
       \edef\@enumctr{enum\romannumeral\the\@enumdepth}\list
       {\csname  label\@enumctr\endcsname}{\setlength{\topsep}{1ex}
\setlength{\itemsep}{0 pt}\usecounter
         {\@enumctr}\def\makelabel##1{\hss\llap{##1}}}\fi}{\endlist}

\renewenvironment{itemize}{\ifnum \@itemdepth >3 \@toodeep\else
\advance\@itemdepth \@ne
\edef\@itemitem{labelitem\romannumeral\the\@itemdepth}
\list{\csname\@itemitem\endcsname}{\setlength{\topsep}{1ex}\setlength
{\itemsep}{0pt}\def\makelabel##1{\hss\llap{##1}}}\fi}
{\endlist}
\def\@seccntformat#1{\csname the#1\endcsname.\quad}

\def\section{\pagebreak[3]\setcounter{prop}{0}\setcounter{equation}{0}\@startsection{section}{1}{\z@}{4ex plus  6ex}{2ex}{\center\reset@font \large\bf}}
\newcommand{\subsect}[1]{\medskip\par\noindent\pagebreak[3]\refstepcounter{subsection}\refstepcounter{prop}{\bf \thesection.\arabic{prop}.\ #1.\ }}
\makeatother

\def\theprop{\thesection.\arabic{prop}}

\renewenvironment{equation}{\refstepcounter{subsection}\refstepcounter
{prop}$$}{\leqno{\bf (\theprop)}$$}

\newenvironment{rem}[1]{\refstepcounter{subsection}\refstepcounter
{prop} \mpn{{\bf \thesection.\arabic{prop}.}\ \ \bf#1~:}}{\smp}
\newenvironment{enonce}[1]{\pagebreak[3]\refstepcounter{prop}\mmpn
{{\bf  \thesection.\arabic{prop}.\ #1.}}\begin{it} }{\end{it}\smp}
\def\thesection{\arabic{section}}

\newcommand{\result}[1]{\begin{enonce}{#1}}
\newcommand{\fresult}{\end{enonce}}
\def\monitem{\medskip\par\noindent$\bullet$ }

\begin{document}
\centerline{\Large\bf The primitive idempotents}\vspace{.2cm}
\centerline{\Large\bf of the $p$-permutation ring}\vspace{.4cm}
\centerline{Serge Bouc and Jacques Th\'evenaz}\vspace{.5cm}
\begin{footnotesize}
{\bf Abstract~:} Let $G$ be a finite group, let $p$ be a prime number, and let $K$ be a field of characteristic 0 and $k$ be a field of characteristic $p$, both large enough. In this note we state explicit formulae for the primitive idempotents of $K\otimes_\Z pp_k(G)$, where $pp_k(G)$ is the ring of $p$-permutation $kG$-modules.\bigskip\par
{\bf AMS Subject Classification~:} 19A22, 20C20.\par
{\bf Key words~:} $p$-permutation, idempotent, trivial source. 
\end{footnotesize}
\section{Introduction}
Let $G$ be a finite group, let $p$ be a prime number, and let $K$ be a field of characteristic 0 and $k$ be a field of characteristic $p$, both large enough. In this note we state explicit formulae for the primitive idempotents of $K\otimes_\Z pp_k(G)$, where $pp_k(G)$ is the ring of $p$-permutation $kG$-modules (also called the trivial source ring).\par
To obtain these formulae, we first use induction and restriction to reduce to the case where $G$ is cyclic modulo $p$, i.e. $G$ has a normal Sylow $p$-subgroup with cyclic quotient. Then we solve the easy and well known case where $G$ is a cyclic $p'$-group. Finally we conclude by using the natural ring homomomorphism from the Burnside ring $B(G)$ of $G$ to $pp_k(G)$, and the classical formulae for the primitive idempotents of $K\otimes_\Z B(G)$.\par
Our formulae are an essential tool in \cite{cartan-mackey}, where Cartan matrices of Mackey algebras are considered, and some invariants of these matrices (determinant, rank) are explicitly computed.
\vspace{-.1cm}
\section{$p$-permutation modules}
\subsect{Notation}
\monitem Throughout the paper, $G$ will be a fixed finite group and $p$ a fixed prime number. We consider a field $k$ of characteristic $p$ and we denote by $kG$ the group algebra of $G$ over~$k$. We assume that $k$ is large enough in the sense that it is a splitting field for every group algebra $k\big(N_G(P)/P\big)$, where $P$ runs through the set of all $p$-subgroups of~$G$.
\monitem We let $K$ be a field of characteristic~0 and we assume that $K$ is large enough in the sense that it contains the values of all the Brauer characters of the groups $N_G(P)/P$, where $P$ runs through the set of all $p$-subgroups of~$G$.\par
We recall quickly how Brauer characters are defined. We let $\overline k$ be an extension of~$k$ containing all the $n$-th roots of unity, where $n$ is the $p'$-part of the exponent of~$G$. We choose an isomorphism $\theta: \mu_n({\overline k}) \to \mu_n(\C)$ from the group of $n$-th roots of unity in~$\overline k$ and the corresponding group in~$\C$.
If $V$ is an $r$-dimensional $kH$-module for the group $H=N_G(P)/P$ and if $s$ is an element of the set $H_{p'}$ of all $p'$-elements of~$H$, the matrix of the action of $s$ on~$V$ has eigenvalues $(\lambda_1, \ldots, \lambda_r)$ in the group $\mu_n({\overline k})$. The {\it Brauer character} $\phi_V$ of $V$ is the central function defined on $H_{p'}$, with values in the field $\Q[\mu_n(\C)]$, sending $s$ to $\sum_{i=1}^r \theta(\lambda_i)$. The actual values of Brauer characters may lie in a subfield of $\Q[\mu_n(\C)]$ and we simply require that $K$ contains all these values.
\begin{rem}{Remark} \label{tensoriel} Let $V$ be as above and let $W$ be a $t$-dimensional $kH$-module. If $s$ has eigenvalues $(\mu_1,\ldots,\mu_t)$ on $W$, its eigenvalues for the diagonal action of $H$ on $V\otimes_kW$ are $(\lambda_i\mu_j)_{1\leq i\leq r,\,1\leq j\leq t}$. It follows that $\phi_{V\otimes_kW}(s)=\sum_{i=1}^r\sum_{j=1}^t\theta(\lambda_i\mu_j)=\phi_V(s)\phi_W(s)$.
\end{rem}
\monitem When $H$ is a subgroup of $G$, and $M$ is a $kG$-module, we denote by $\Res_H^GM$ the $kH$-module obtained by restricting the action of $G$ to $H$. When $L$ is a $kH$-module, we denote by $\Ind_H^GL$ the induced $kG$-module. 
\monitem When $M$ is a $kG$-module, and $P$ is a subgroup of $G$, the $k$-vector space of fixed points of $P$ on $M$ is denoted by $M^P$. When $Q\leq P$ are subgroups of $G$, the {\em relative trace} is the map $\tr_Q^P:M^Q\to M^P$ defined by $\tr_Q^P(m)=\sum_{x\in[P/Q]}x\cdot m$. 
\monitem When $M$ is a $kG$-module, the {\em Brauer quotient} of $M$ at $P$ is the $k$-vector space
$$M[P]=M^P/\sum_{Q<P}\tr_Q^PM^Q\mpoint$$
This $k$-vector space has a natural structure of $k\sur{N}_G(P)$-module, where as usual $\sur{N}_G(P)=N_G(P)/P$. It is equal to zero if $P$ is not a $p$-group.
\monitem If $P$ is a normal $p$-subgroup of $G$ and $M$ is a $k(G/P)$-module, denote by $\Inf_{G/P}^GM$ the $kG$-module obtained from $M$ by inflation to $G$. Then there is an isomorphism
$$(\Inf_{G/P}^GM)[P]\cong M$$
of $k(G/P)$-modules.
\monitem When $G$ acts on a set $X$ (on the left), and $x,y\in X$, we write $x=_Gy$ if $x$~and $y$ are in the same $G$-orbit. We denote by $[G\dom X]$ a set of representatives of $G$-orbits on $X$, and by $X^G$ the set of fixed points of $G$ on $X$. For $x\in X$, we denote by $G_x$ its stabilizer in $G$.
\subsect{Review of $p$-permutation modules}
We begin by recalling some definitions and basic results. We refer to~\cite{brouepperm}, and to \cite{benson1} Sections 3.11 and 5.5 for details~:
\result{Definition} A {\em permutation} $kG$-module is a $kG$-module admitting a $G$-invariant $k$-basis. A {\em $p$-permutation} $kG$-module $M$ is a $kG$-module such that $\Res_S^GM$ is a permutation $kS$-module, where $S$ is a Sylow $p$-subgroup of~$G$.
\fresult
The $p$-permutation $kG$-modules are also called {\em trivial source modules}, because the indecomposable ones coincide with the indecomposable modules having a trivial source (see \cite{brouepperm} 0.4). Moreover, the $p$-permutation modules also coincide with the direct summands of permutation modules (see \cite{benson1}, Lemma 3.11.2).
\result{Proposition}\begin{enumerate}
\item If $H$ is a subgroup of $G$, and $M$ is a $p$-permutation $kG$-module, then the restriction $\Res_H^GM$ of $M$ to $H$ is a $p$-permutation $kH$-module.
\item If $H$ is a subgroup of $G$, and $L$ is a $p$-permutation $kH$-module, then the induced module $\Ind_H^GL$ is a $p$-permutation $kG$-module.
\item If $N$ is a normal subgroup of $G$, and $L$ is a $p$-permutation $k(G/N)$-module, the inflated module $\Inf_{G/N}^GL$ is a $p$-permutation $kG$-module.
\item If $P$ is a $p$-group, and $M$ is a permutation $kP$-module with $P$-invariant basis $X$, then the image of the set $X^P$ in $M[P]$ is a $k$-basis of $M[P]$.
\item If $P$ is a $p$-subgroup of $G$, and $M$ is a $p$-permutation $kG$-module, then the Brauer quotient $M[P]$ is a $p$-permutation $k\sur{N}_G(P)$-module.
\item If $M$ and $N$ are $p$-permutation $kG$-modules, then their tensor product $M\otimes_kN$ is again a $p$-permutation $kG$-module.
\end{enumerate}
\fresult
\pf Assertions 1,2,3, and 6 are straightforward consequences of the same assertions for {\em permutation} modules. For Assertion~4, see \cite{brouepperm} 1.1.(3). Assertion~5 follows easily from Assertion~4 (see also \cite{brouepperm} 3.1).\endpf
This leads to the following definition~:
\result{Definition} The {\em $p$-permutation ring} $pp_k(G)$ is the Grothendieck group of the category of $p$-permutation $kG$-modules, with relations corresponding to direct sum decompositions, i.e. $[M]+[N]=[M\oplus N]$. The ring structure on $pp_k(G)$ is induced by the tensor product of modules over $k$. The identity element of $pp_k(G)$ is the class of the trivial $kG$-module $k$.
\fresult
As the Krull-Schmidt theorem holds for $kG$-modules, the additive group $pp_k(G)$ is a free (abelian) group on the set of isomorphism classes of indecomposable $p$-permutation $kG$-modules. These modules have the following properties~:
\result{Theorem} {\rm[ \cite{brouepperm} Theorem 3.2]} \label{broue}\begin{enumerate} 
\item The vertices of an indecomposable $p$-permutation $kG$-module $M$ are the maximal $p$-subgroups $P$ of $G$ such that $M[P]\neq\zero$.
\item An indecomposable $p$-permutation $kG$-module has vertex $P$ if and only if $M[P]$ is a non-zero projective $k\sur{N}_G(P)$-module.
\item The correspondence $M\mapsto M[P]$ induces a bijection between the isomorphism classes of indecomposable $p$-permutation $kG$-modules with vertex~$P$ and the isomorphism classes of indecomposable projective $k\sur{N}_G(P)$-modules.
\end{enumerate}
\fresult
\result{Notation} Let $\mathcal{P}_{G,p}$ denote the set of pairs $(P,E)$, where $P$ is a $p$-subgroup of $G$, and $E$ is an indecomposable projective $k\sur{N}_G(P)$-module. The group $G$ acts on $\mathcal{P}_{G,p}$ by conjugation, and we denote by $[\mathcal{P}_{G,p}]$ a set of representatives of $G$-orbits on $\mathcal{P}_{G,p}$.\par
For $(P,E)\in\mathcal{P}_{G,p}$, let $M_{P,E}$ denote the (unique up to isomorphism) indecomposable $p$-permutation $kG$-module such that $M_{P,E}[P]\cong E$. 
\fresult
\result{Corollary} \label{base}The classes of the modules $M_{P,E}$, for $(P,E)\in[\mathcal{P}_{G,p}]$ form a $\Z$-basis of $pp_k(G)$.
\fresult
\result{Notation} The operations $\Res_H^G$, $\Ind_H^G$, $\Inf_{G/N}^G$ extend linearly to maps between the corresponding $p$-permutations rings, denoted with the same symbol. 
\fresult
The maps $\Res_H^G$ and $\Inf_{G/N}^G$ are ring homomorphisms, whereas $\Ind_H^G$ is not in general. Similarly~:
\result{Proposition} \label{ring hom}Let $P$ be a $p$-subgroup of $G$. Then the correspondence $M\mapsto M[P]$ induces a ring homomorphism $\Br_P^G:pp_k(G)\to pp_k\big(\sur{N}_G(P)\big)$.
\fresult
\pf Let $M$ and $N$ be $p$-permutation $kG$-modules. The canonical bilinear map $M\times N\to M\otimes_kN$ is $G$-equivariant, hence it induces a bilinear map $\beta_P: M[P]\times N[P]\to (M\otimes_k N)[P]$ (see \cite{brouepperm} 1.2), which is $\sur{N}_G(P)$-equivariant. Now if $X$ is a $P$-invariant $k$-basis of $M$, and $Y$ a $P$-invariant $k$-basis of $N$, then $X\times Y$ is a $P$-invariant basis of $M\otimes_k N$. The images of the sets $X^P$, $Y^P$, and $(X\times Y)^P$ are bases of $M[P]$, $N[P]$, and $(M\otimes_k N)[P]$, respectively, and the restriction of $\beta_P$ to theses bases is the canonical bijection $X^P\times Y^P\to (X\times Y)^P$. It follows that $\beta_P$ induces an isomorphism $M[P]\otimes_kN[P]\to (M\otimes_k N)[P]$ of $k\sur{N}_G(P)$-modules. Proposition~\ref{ring hom} follows.\endpf 
\result{Notation}
 Let $\mathcal{Q}_{G,p}$ denote the set of pairs $(P,s)$, where $P$ is a $p$-subgroup of $G$, and $s$ is a $p'$-element of $\sur{N}_G(P)$. The group $G$ acts on $\mathcal{Q}_{G,p}$, and we denote by $[\mathcal{Q}_{G,p}]$ a set of representatives of $G$-orbits on $\mathcal{Q}_{G,p}$. \par
If $(P,s)\in \mathcal{Q}_{G,p}$, we denote by $N_G(P,s)$ the stabilizer of $(P,s)$ in $G$, and by $\gen{Ps}$ the subgroup of $N_G(P)$ generated by $Ps$ (i.e. the inverse image in $N_G(P)$ of the cyclic group $\gen{s}$ of $\sur{N}_G(P)$).\par
\fresult
\pagebreak[3]
\begin{rem}{Remarks} \begin{itemize}
\item When $H$ is a subgroup of $G$, there is a natural inclusion of $\mathcal{Q}_{H,p}$ into $\mathcal{Q}_{G,p}$, as $\sur{N}_H(P)\leq \sur{N}_G(P)$ for any $p$-subgroup $P$ of $H$. We will consider $\mathcal{Q}_{H,p}$ as a subset of $\mathcal{Q}_{G,p}$.
\item When $(P,s)\in \mathcal{Q}_{G,p}$, the group $N_G(P,s)$ is the set of elements $g$ in $N_G(P)$ whose image in $\sur{N}_G(P)$ centralizes $s$. In other words, there is a short exact sequence of groups
\begin{equation}
\label{ses}\un\to P\to N_G(P,s)\to C_{\sur{N}_G(P)}(s)\to \un\mpoint
\end{equation}
In particular $N_G(P,s)$ is a subgroup of $N_G(\gen{Ps})$.
\end{itemize}
\end{rem}
\result{Notation} Let $(P,s)\in\mathcal{Q}_{G,p}$. Let $\tau^G_{P,s}$ denote the additive map from $pp_k(G)$ to $K$ sending the class of a $p$-permutation $kG$-module $M$ to the value at $s$ of the Brauer character of the $\sur{N}_G(P)$-module $M[P]$.
\fresult
\begin{rem}{Remarks} \begin{itemize}
\item \label{tau Res}It is clear that $\tau_{P,s}^G(M)$ only depends on the restriction of~$M$ to the group ${<}Ps{>}$. In other words
$$\tau_{P,s}^G=\tau_{P,s}^{{<}Ps{>}}\circ\Res_{{<}Ps{>}}^G\mpoint$$
Furthermore, it is clear from the definition that
\begin{equation} \label{factor tau}\tau_{P,s}^G=\tau_{\un,s}^{\gen{Ps}/P}\circ\Br_P^{\gen{Ps}}\circ\Res_{\gen{Ps}}^G\mpoint
\end{equation}
\item It is easy to see that $\tau_{P,s}^G$ only depends on the $G$-orbit of $(P,s)$, that is, $\tau_{P^g,s^g}^G = \tau_{P,s}^G$ for every $g\in G$.
\end{itemize}
\end{rem}
The following proposition is Corollary 5.5.5 in \cite{benson1}, but our construction of the species is slightly different (but equivalent, of course). For this reason, we sketch an independent proof~:
\result{Proposition} \begin{enumerate}
\item The map $\tau^G_{P,s}$ is a ring homomorphism $pp_k(G)\to K$ and extends to a $K$-algebra homomorphism (a species) $\tau_{P,s}^G : K\otimes_\Z pp_k(G) \to K$.
\item The set $\{\tau^G_{P,s}\mid (P,s)\in[\mathcal{Q}_{G,p}]\}$ is the set of all distinct species from $K\otimes_\Z pp_k(G)$ to $K$. These species induce a $K$-algebra isomorphism
$$T=\prod_{(P,s)\in[\mathcal{Q}_{G,p}]}\tau^G_{P,s}:K\otimes_\Z pp_k(G)\to\prod_{(P,s)\in[\mathcal{Q}_{G,p}]}K\mpoint$$
\end{enumerate}
\vspace{-.8cm}
\fresult
\pf By \ref{factor tau}, to prove Assertion 1, it suffices to prove that $\tau_{\un,s}^{\gen{Ps}/P}$ is a ring homomorphism, since both $\Res_{\gen{Ps}}^G$ and $\Br_P^{\gen{Ps}}$ are ring homomorphisms. In other words, we can assume that $P=\un$. Now the value of $\tau_{\un,s}^G$ on the class of a $kG$-module $M$ is the value $\phi_M(s)$ of the Brauer character of $M$ at~$s$, so Assertion 1 follows from Remark~\ref{tensoriel}. \par
For Assertion 2, it suffices to prove that $T$ is an isomorphism. Since $[\mathcal{P}_{G,p}|$ and $[\mathcal{Q}_{G,p}]$ have the same cardinality, the matrix $\mathcal{M}$ of $T$ is a square matrix. Let $(P,E)\in\mathcal{P}_{G,p}$, and $(Q,s)\in\mathcal{Q}_{G,p}$. Then $\tau_{Q,s}(M_{P,E})$ is equal to zero if $Q$ is not contained in $P$ up to $G$-conjugation, because in this case $M_{P,E}(Q)=\zero$ by Theorem~\ref{broue}. It follows that $\mathcal{M}$ is block triangular. As moreover $M_{P,E}[P]\cong E$, we have that $\tau_{P,s}(M_{P,E})=\phi_E(s)$. This means that the diagonal block of $\mathcal{M}$
corresponding to $P$ is the matrix of Brauer characters of projective $k\sur{N}_G(P)$-modules, and these are linearly independent by Lemma 5.3.1 of~\cite{benson1}. It follows that all the diagonal blocks of $\mathcal{M}$ are non singular, so $\mathcal{M}$ is invertible, and $T$ is an isomorphism.
\endpf
\vspace{-.2cm}
\result{Corollary} The algebra $K\otimes_\Z pp_k(G)$ is a split semisimple commutative $K$-algebra. Its primitive idempotents $F_{P,s}^G$ are indexed by $[\mathcal{Q}_{G,p}]$, and the idempotent $F_{P,s}^G$ is characterized by 
$$\forall (R,u)\in \mathcal{Q}_{G,p},\;\;\tau^G_{R,u}(F_{P,s}^G)=\left\{\begin{array}{cl}1&\hbox{if}\;(R,u)=_G(P,s)\\0&\hbox{otherwise.}\end{array}\right. $$
\fresult
\vspace{-.8cm}
\section{Restriction and induction}\nopagebreak
\vspace{-.2cm}
\result{Proposition} \label{res idemp}Let $H\leq G$, and $(P,s)\in\mathcal{Q}_{G,p}$. Then
$$\Res_H^GF_{P,s}^G=\sum_{(Q,t)}F^H_{Q,t}\mvirg$$
where $(Q,t)$ runs through a set of representatives of $H$-conjugacy classes of $G$-conjugates of $(P,s)$ contained in $H$.
\fresult
\pf Indeed, as $\Res_H^G$ is an algebra homomorphism, the element $\Res_H^GF_{P,s}^G$ is an idempotent of $K\otimes_\Z pp_k(H)$, hence it is equal to a sum of some distinct primitive idempotents $F_{Q,t}^H$. The idempotent $F_{Q,t}^H$ appears in this decomposition if and only if $\tau_{Q,t}^H(\Res_H^GF_{P,s}^G)=1$. By Remark~\ref{tau Res}
\begin{eqnarray*}
\tau_{Q,t}^H(\Res_H^GF_{P,s}^G)&=&\tau_{Q,t}^{{<}Qt{>}}(\Res_{\gen{Qt}}^H\Res_H^GF_{P,s}^G)\\
&=&\tau_{Q,t}^{\gen{Qt}}(\Res_{\gen{Qt}}^GF_{P,s}^G)\\
&=&\tau_{Q,t}^G(F_{P,s}^G)\mpoint
\end{eqnarray*}
Now $\tau_{Q,t}^G(F_{P,s}^G)$ is equal to 1 if and only if $(Q,t)$ and $(P,s)$ are $G$-conjugate. This completes the proof.\endpf
\result{Proposition} \label{ind idemp}Let $H\leq G$, and $(Q,t)\in\mathcal{Q}_{H,p}$. Then
$$\Ind_H^GF_{Q,t}^H=|N_G(Q,t):N_H(Q,t)|F_{Q,t}^G\mpoint$$
\fresult
\pf Since $K\otimes_\Z pp_k(G)$ is a split semisimple commutative $K$-algebra, any element $X$ in $K\otimes_\Z pp_k(G)$ can be written
\begin{equation}\label{decomp}
X=\sum_{(P,s)\in[\mathcal{Q}_{G,p}]}\tau_{P,s}^G(X)F_{P,s}^G\mvirg
\end{equation}
and moreover for any $(P,s)\in\mathcal{Q}_{G,p}$
$$\tau_{P,s}^G(X)F_{P,s}^G=X\cdot F_{P,s}^G\mpoint$$
Setting $X=\Ind_H^GF_{Q,t}^H$ in this equation gives
\begin{eqnarray*}
\tau_{P,s}^G(\Ind_H^GF_{Q,t}^H)F_{P,s}^G&=&(\Ind_H^GF_{Q,t}^H)\cdot F_{P,s}^G\\
&=&\Ind_H^G(F_{Q,t}^H\cdot\Res_H^GF_{P,s}^G)\mpoint
\end{eqnarray*}
By Proposition~\ref{res idemp}, the element $\Res_H^GF_{P,s}^G$ is equal to the sum of the distinct idempotents $F_{P^y,s^y}^H$ associated to elements $y$ of $G$ such that $\gen{Ps}^y\leq H$. The product $F_{Q,t}^H\cdot F_{P^y,s^y}^H$ is equal to zero, unless $(Q,t)$ is $H$-conjugate to $(P^y,s^y)$, which implies that $(Q,t)$ and $(P,s)$ are $G$-conjugate. It follows that the only non zero term in the right hand side of Equation~\ref{decomp} is the term corresponding to $(Q,t)$. Hence
$$\Ind_H^GF_{Q,t}^H=\tau_{Q,t}^G(\Ind_H^GF_{Q,t}^H)F_{Q,t}^G\mpoint$$
Now by Remark~\ref{tau Res} and the Mackey formula
\begin{eqnarray*}
\tau_{Q,t}^G(\Ind_H^GF_{Q,t}^H)&=&\tau_{Q,t}^{\gen{Qt}}(\Res_{\gen{Qt}}^G\Ind_H^GF_{Q,t}^H)\\
&=&\tau_{Q,t}^{\gen{Qt}}\big(\sum_{x\in\gen{Qt}\dom G/H}\Ind_{\gen{Qt}\cap{^xH}}^{\gen{Qt}}{^x\Res}_{\gen{Qt}^x\cap H}^HF_{Q,t}^H\big)\mpoint\\
\end{eqnarray*}
By Proposition~\ref{res idemp}, the element $\Res_{\gen{Qt}^x\cap H}^HF_{Q,t}^H$ is equal to the sum of the distinct idempotents $F_{Q^y,t^y}^{\gen{Qt}^x\cap H}$ corresponding to elements $y\in H$ such that $\gen{Qt}^y\leq \gen{Qt}^x\cap H$. This implies $\gen{Qt}^y=\gen{Qt}^x$, i.e. $y\in N_G(\gen{Qt})x$, thus $x\in N_G(\gen{Qt})\cdot H$. This gives
\begin{eqnarray*}
\tau_{Q,t}^G(\Ind_H^GF_{Q,t}^H)&=&\tau_{Q,t}^{\gen{Qt}}\big(\sumb{x\in N_G(\gen{Qt})H/H}{\rule{0ex}{1.2ex}y\in N_H(Q,t)\dom N_G(\gen{Qt})x}{^xF}_{Q^y,t^y}^{\gen{Qt}}\big)\\
&=&\sumb{x\in N_G(\gen{Qt})/N_H(\gen{Qt})}{\rule{0ex}{1.2ex}y\in N_H(Q,t)\dom N_G(\gen{Qt})x}\tau_{Q,t}^{\gen{Qt}}(F_{Q^{yx^{-1}},t^{yx^{-1}}}^{\gen{Qt}})\\
&=&\sum_{z\in N_H(Q,t)\dom N_G(\gen{Qt})}\tau_{Q,t}^{\gen{Qt}}(F_{Q^z,t^z}^{\gen{Qt}})\mvirg
\end{eqnarray*}
where $z=yx^{-1}$. Finally $\tau_{Q,t}^{\gen{Qt}}(F_{Q^z,t^z}^{\gen{Qt}})$ is equal to 1 if $(Q^z,t^z)$ is conjugate to $(Q,t)$ in $\gen{Qt}$, and to zero otherwise. \par
If $u\in \gen{Qt}$ is such that $(Q^z,t^z)^u=(Q,t)$, then $zu\in N_G(Q,t)$. But since $[\gen{Qt},t]\leq Q$, we have $\gen{Qt}\leq N_G(Q,t)$, so $u\in N_G(Q,t)$, hence $z\in N_G(Q,t)$, and $(Q^z,t^z)=(Q,t)$. It follows that 
$$\tau_{Q,t}^G(\Ind_H^GF_{Q,t}^H)=|N_G(Q,t):N_H(Q,t)|\mvirg$$
which completes the proof of the proposition.\endpf
\result{Corollary} \label{reduction} Let $(P,s)\in\mathcal{Q}_{G,p}$. Then
$$F_{P,s}^G=\frac{|s|}{|C_{\sur{N}_G(P)}(s)|}\Ind_{\gen{Ps}}^GF_{P,s}^{\gen{Ps}}\mpoint$$
\fresult
\pf Apply Proposition~\ref{ind idemp} with $(Q,t)=(P,s)$ and $H=\gen{Ps}$. Then $N_H(Q,t)=\gen{Ps}$, thus by Exact sequence~\ref{ses}
$$|N_G(Q,t):N_H(Q,t)|=\frac{|P||C_{\sur{N}_G(P)}(s)|}{|P||s|}=\frac{|C_{\sur{N}_G(P)}(s)|}{|s|}\mvirg$$
and the corollary follows.\endpf
\section{Idempotents}
It follows from Corollary \ref{reduction} that, in order to find formulae for the primitive idempotents $F_{P,s}^G$ of $K\otimes_\Z pp_k(G)$, it suffices to find the formula expressing the idempotent $F_{P,s}^{\gen{Ps}}$. In other words, we can assume that $G=\gen{Ps}$, i.e. that $G$ has a normal Sylow $p$-subgroup $P$ with cyclic quotient generated by~$s$.
\subsect{A morphism from the Burnside ring}
When $G$ is an arbitrary finite group, there is an obvious ring homomorphism $\mathcal{L}_G$ from the Burnside ring $B(G)$ to $pp_k(G)$, induced by the {\em linearization} operation, sending a finite $G$-set $X$ to the permutation module $kX$, which is obviously a $p$-permutation module. This morphism also commutes with restriction and induction~: if $H\leq G$, then
\begin{equation}
\label{commute}\mathcal{L}_H\circ\Res_H^G=\Res_H^G\circ\mathcal{L}_G\mvirg\hspace{1cm}\mathcal{L}_G\circ\Ind_H^G=\Ind_H^G\circ\mathcal{L}_H\mpoint\end{equation}
Indeed, for any $G$-set $X$, the $kH$-modules $k\Res_H^GX$ and $\Res_H^G(kX)$ are isomorphic, and for any $H$-set $Y$, the $kG$-modules $k\Ind_H^GY$ and $\Ind_H^G(kY)$ are isomorphic.\par
Similarly, when $P$ is a $p$-subgroup of $G$, the ring homomorphism $\Phi_P:B(G)\to B\big(\sur{N}_G(P)\big)$ induced by the operation $X\mapsto X^P$ on $G$-sets, is compatible with the Brauer morphism $\Br_P^G:pp_k(G)\to pp_k\big(\sur{N}_G(P)\big)$~:
\begin{equation}\label{commute Brauer}
\mathcal{L}_{\sur{N}_G(P)}\circ \Phi_P=\Br_P^G\circ \mathcal{L}_G\mpoint
\end{equation}
This is because for any $G$-set $X$, the $k\sur{N}_G(P)$-modules $k(X^P)$ and $(kX)[P]$ are isomorphic.\par
The ring homomorphism $\mathcal{L}_G$ extends linearly to a $K$-algebra homomorphism $K\otimes_\Z B(G)\to K\otimes_\Z pp_k(G)$, still denoted by $\mathcal{L}_G$. The algebra $K\otimes_\Z B(G)$ is also a split semisimple commutative $K$-algebra. Its species are the $K$-algebra maps
$$K\otimes_\Z B(G)\to K,\;\;X\mapsto|X^H|\mvirg$$
where $H$ runs through the set of all subgroups of $G$ up to conjugation. Here we denote by $|X^H|$ the number of $H$-fixed points of a $G$-set~$X$ and this notation is then extended $K$-linearly to any $X\in K\otimes_\Z B(G)$. The primitive idempotents $e_H^G$ of $K\otimes_\Z B(G)$ are indexed by subgroups $H$ of $G$, up to conjugation. They are given by the following formulae, found by Gluck (\cite{gluck}) and later independently by Yoshida~(\cite{yoshidaidemp})~:
\begin{equation}\label{gluck yoshida}
e_H^G=\frac{1}{|N_G(H)|}\sum_{L\leq H}|L|\,\mu(L,H)\,G/L\mvirg
\end{equation}
where $\mu$ denotes the M\"obius function of the poset of subgroups of $G$. The idempotent $e_H^G$ is characterized by the fact that for any $X\in K\otimes_\Z B(G)$
$$X\cdot e_H^G=|X^H|e_H^G\mpoint$$
\begin{rem}{Remark} \label{idemp induit}
Since $|X^H|$ only depends on $\Res_H^GX$, it follows in particular that $X$ is a scalar multiple of the ``top" idempotent $e_G^G$ if and only if $\Res_H^GX=\nolinebreak 0$ for any proper subgroup $H$ of $G$. In particular, if $N$ is a normal subgroup of $G$, then
\begin{equation}\label{points fixes}
(e_G^G)^N=e_{G/N}^{G/N}\mpoint
\end{equation}
This is because for any proper subgroup $H/N$ of $G/N$
$$\Res_{H/N}^{G/N}(e_G^G)^N=(\Res_H^Ge_G^G)^N=0\mpoint$$
So $(e_G^G)^N$ is a scalar multiple of $e_{G/N}^{G/N}$. As it is also an idempotent, it is equal to 0 or $e_{G/N}^{G/N}$. Finally
$$|\big((e_G^G)^N\big)^{G/N}|=|(e_G^G)^G|=1\mvirg$$
so $(e_G^G)^N$ is non zero.
\end{rem}
\subsect{The case of a cyclic $p'$-group}
Suppose that $G$ is a cyclic $p'$-group, of order $n$, generated by an element $s$. In this case, there are exactly $n$ group homomorphisms from $G$ to the multiplicative group $k^\times$ of $k$. For each of these group homomorphisms~$\varphi$, let $k_\varphi$ denote the $kG$-module $k$ on which the generator $s$ acts by multiplication by $\varphi(s)$. As $G$ is a $p'$-group, this module is simple and projective. The (classes of the) modules $k_\varphi$, for $\varphi\in \widehat{G}=\Hom(G,k^\times)$, form a basis of $pp_k(G)$. \par
Since moreover for $\varphi, \psi\in \widehat{G}$, the modules $k_\varphi\otimes_kk_\psi$ and $k_{\varphi\psi}$ are isomorphic, the algebra $K\otimes_\Z pp_k(G)$ is isomorphic to the group algebra of the group $\widehat{G}$. This leads to the following classical formula~:
\result{Lemma} \label{cyclique}Let $G$ be a cyclic $p'$-group. Then for any $t\in G$,
$$F_{\un,t}^G=\frac{1}{n}\sum_{\varphi\in \widehat{G}}\tilde{\varphi}(t^{-1})k_\varphi\mvirg
$$
where $\tilde{\varphi}$ is the Brauer character of $k_\varphi$. 
\fresult
\pf Indeed for $s,t\in G$
$$\tau_{\un,t}^G\big(\frac{1}{n}\sum_{\varphi\in \widehat{G}}\tilde{\varphi}(s^{-1})k_\varphi\big)=\frac{1}{n}\sum_{\varphi\in \widehat{G}}\tilde{\varphi}(s^{-1})\tilde{\varphi}(t)=\delta_{s,t}\mvirg$$
where $\delta_{s,t}$ is the Kronecker symbol.\endpf
\subsect{The case $G=\gen{Ps}$}
Suppose now more generally that $G=\gen{Ps}$, where $P$ is a normal Sylow $p$-subgroup of $G$ and $s$ is a $p'$-element. In this case, by Proposition~\ref{res idemp}, the restriction of $F_{P,s}^G$ to any proper subgroup of~$G$ is equal to zero. Moreover, since $N_G(P,t)=G$ for any $t\in G/P$, the conjugacy class of the pair $(P,t)$ reduces to $\{(P,t)\}$.
\result{Lemma} \label{cyclique mod p}Suppose $G=\gen{Ps}$, and set $E_G^G=\mathcal{L}_G(e_G^G)$. Then
$$E_G^G=\sum_{\gen{t}=\gen{s}}F_{P,t}^G\mpoint$$
\fresult
\pf By \ref{commute} and by Remark~\ref{idemp induit}, the restriction of $E_G^G$ to any proper subgroup of $G$ is equal to zero. Let $(Q,t)\in\mathcal{Q}_{G,p}$, such that the group $L=\gen{Qt}$ is a proper subgroup of $G$. By Proposition~\ref{ind idemp}, there is a rational number $r$ such that
$$F_{Q,t}^G=r\,\Ind_L^GF_{Q,t}^L\mpoint$$
It follows that
$$E_G^G\cdot F_{Q,t}^G=r\,\Ind_L^G\big((\Res_L^GE_G^G)\cdot F_{Q,t}^L\big)=0\mpoint$$
Now $E_G^G$ is an idempotent of $K\otimes_\Z pp_k(G)$, hence is it a sum of some of the primitive idempotents $F_{Q,t}^G$ associated to pairs $(Q,t)$ for which $\gen{Qt}=G$. This condition is equivalent to $Q=P$ and $\gen{t}=\gen{s}$.\par
It remains to show that all these idempotents $F_{P,t}^G$ appear in the decomposition of $E_G^G$, i.e. equivalently that $\tau_{P,t}^G(E_G^G)=1$ for any generator $t$ of~$\gen{s}$. Now by \ref{points fixes} and Remark~\ref{tau Res}
$$\tau_{P,t}^G(E_G^G)=\tau_{\un,t}^{G/P}\big(\Br_P^G(E_G^G)\big)=\tau_{\un,t}^{\gen{s}}(E_{\gen{s}}^{\gen{s}})\mpoint$$
Now the value at $t$ of the Brauer character of a permutation module $kX$ is equal to the number of fixed points of $t$ on $X$. By $K$-linearity, this gives
$$\tau_{\un,t}^{\gen{s}}(E_{\gen{s}}^{\gen{s}})=|(e_{\gen{s}}^{\gen{s}})^t|\mvirg$$
and this is equal to 1 if $t$ generates $\gen{s}$, and to 0 otherwise, as was to be shown.\endpf
\result{Proposition} \label{cas G=Ps}Let $(P,s)\in \mathcal{Q}_{G,p}$, and suppose that $G=\gen{Ps}$. Then
$$F_{P,s}^G=E_G^G\cdot\Inf_{G/P}^GF_{\un,s}^{G/P}\mpoint$$
\fresult
\pf Set $E_s=E_G^G\cdot\Inf_{G/P}^GF_{\un,s}^{G/P}$. Then $E_s$ is an idempotent of $K\otimes_\Z pp_k(G)$, as it is the product of two (commuting) idempotents. Let $(Q,t)\in\mathcal{Q}_{G,p}$. If $\gen{Qt}\neq G$, then $\tau_{Q,t}^G(E_G^G)=0$ by Lemma~\ref{cyclique mod p}, thus $\tau_{Q,t}^G(E_s)=0$. And if $\gen{Qt}=G$, then $Q=P$ and $\gen{t}=\gen{s}$. In this case
$$\tau_{Q,t}^G(E_s)=\tau_{P,t}^G(E_G^G)\cdot\tau_{P,t}^G(\Inf_{G/P}^GF_{\un,s}^{G/P})\mpoint$$
By Lemma~\ref{cyclique mod p}, the first factor in the right hand side is equal to 1. The second factor is equal to
\begin{eqnarray*}
\tau_{P,t}^G(\Inf_{G/P}^GF_{\un,s}^{G/P})&=&\tau_{\un,t}^{G/P}\Br_P^G(\Inf_{G/P}^GF_{\un,s}^{G/P})\\
&=&\tau_{\un,t}^{G/P}(F_{\un,s}^{G/P})=\delta_{t,s}\mvirg
\end{eqnarray*}
where $\delta_{t,s}$ is the Kronecker symbol. Hence $\tau_{P,t}^G(E_s)=\delta_{t,s}$, and this completes the proof.\endpf
\result{Theorem} Let $G$ be a finite group, and let $(P,s)\in \mathcal{Q}_{G,p}$. Then the primitive idempotent $F_{P,s}^G$ of the $p$-permutation algebra $K\otimes_\Z pp_k(G)$ is given by the following formula~:
$$F_{P,s}^G=\frac{1}{|P||s||C_{\sur{N}_G(P)}(s)|}\sumc{\varphi\in\widehat{\gen{s}}}{\rule{0ex}{1.3ex}L\leq \gen{Ps}}{\rule{0ex}{1.3ex}PL=\gen{Ps}}\tilde{\varphi}(s^{-1})|L|\mu(L,\gen{Ps})\,\Ind_L^Gk_{L,\varphi}^{\gen{Ps}}\mvirg$$
where $k_{L,\varphi}^{\gen{Ps}}=\Res_L^{\gen{Ps}}\Inf_{\gen{s}}^{\gen{Ps}}k_\varphi$.
\fresult
\pf By Corollary~\ref{reduction}, and Proposition~\ref{cas G=Ps}
$$F_{P,s}^G=\frac{|s|}{|C_{\sur{N}_G(P)}(s)|}\Ind_{\gen{Ps}}^G(E_{\gen{Ps}}^{\gen{Ps}}\cdot\Inf_{\gen{s}}^{\gen{Ps}}F_{\un,s}^{\gen{s}})\mpoint$$
By Equation~\ref{gluck yoshida}, this gives
$$F_{P,s}^G\!=\!\frac{|s|}{|C_{\sur{N}_G(P)}(s)|}\Ind_{\gen{Ps}}^G\frac{1}{|P||s|}\!\!\sum_{L\leq \gen{Ps}}\!\!|L|\mu(L,\gen{Ps})\Ind_L^{\gen{Ps}}k\cdot\Inf_{\gen{s}}^{\gen{Ps}}F_{\un,s}^{\gen{s}}\mpoint$$
Moreover for each $L\leq\gen{Ps}$
\begin{eqnarray*}
\Ind_L^{\gen{Ps}}k\cdot\Inf_{\gen{s}}^{\gen{Ps}}F_{\un,s}^{\gen{s}}&\cong&\Ind_L^{\gen{Ps}}(\Res_L^{\gen{Ps}}\Inf_{\gen{s}}^{\gen{Ps}}F_{\un,s}^{\gen{s}})\\
&\cong&\Ind_L^{\gen{Ps}}\Inf_{L/L\cap P}^L\Iso_{LP/P}^{L/L\cap P}\Res_{LP/P}^{\gen{s}}F_{\un,s}^{\gen{s}}\mpoint
\end{eqnarray*}
Here we have used the fact that if $L$ and $P$ are subgroups of a group $H$, with $P\normal H$, then there is an isomorphism of functors
$$\Res_L^H\circ\Inf_{H/P}^H\cong \Inf_{L/L\cap P}^L\circ\Iso_{LP/P}^{L/L\cap P}\circ\Res_{LP/P}^{H/P}\mvirg$$
which follows from the isomorphism of $(L,H/P)$-bisets
$$H\times_H(H/P)\cong {_L(}H/P)_{H/P}\cong (L/L\cap P)\times_{L/L\cap P}(LP/P)\times_{LP/P}(H/P)\mpoint$$
Now Proposition~\ref{res idemp} implies that $\Res_{LP/P}^{\gen{s}}F_{\un,s}^{\gen{s}}=0$  if $LP/P\neq \gen{s}$, i.e. equivalently if $PL\neq\gen{Ps}$. It follows that
$$
F_{P,s}^G\!=\!\frac{1}{|P||C_{\sur{N}_G(P)}(s)|}\!\sumb{L\leq \gen{Ps}}{\rule{0ex}{1.3ex}PL=\gen{Ps}}\!\!\!|L|\mu(L,\gen{Ps})\Ind_L^{G}(\Res_L^{\gen{Ps}}\Inf_{\gen{s}}^{\gen{Ps}}F_{\un,s}^{\gen{s}})\mpoint
$$
By Lemma~\ref{cyclique}, this gives
$$
F_{P,s}^G\!=\!\frac{1}{|P||s||C_{\sur{N}_G(P)}(s)|}\!\sumc{\varphi\in\widehat{\gen{s}}}{\rule{0ex}{1.3ex}L\leq \gen{Ps}}{\rule{0ex}{1.3ex}PL=\gen{Ps}}\!\!\!\tilde{\varphi}(s^{-1})|L|\mu(L,\gen{Ps})\Ind_L^{G}k_{L,\varphi}^{\gen{Ps}}\mvirg
$$
where $k_{L,\varphi}^{\gen{Ps}}=\Res_L^{\gen{Ps}}\Inf_{\gen{s}}^{\gen{Ps}}k_\varphi$, as was to be shown.\endpf

\bigskip
\noindent
Serge Bouc, CNRS-LAMFA, Universit\'e de Picardie - Jules Verne,\\
33, rue St Leu, F-80039 Amiens Cedex~1, France.\\
{\tt serge.bouc@u-picardie.fr}

\medskip
\noindent
Jacques Th\'evenaz, Institut de G\'eom\'etrie, Alg\`ebre et Topologie,\\
EPFL, B\^atiment BCH, CH-1015 Lausanne, Switzerland.\\
{\tt Jacques.Thevenaz@epfl.ch}

\end{document}